# Surrogate-based toll optimization in a large-scale heterogeneously congested network

Ziyuan Gu, S. Travis Waller, Meead Saberi*

*School of Civil and Environmental Engineering, Research Center for Integrated Transport Innovation (rCITI), University of New South Wales, Sydney NSW 2052, Australia*

**Abstract:** Toll optimization in a large-scale dynamic traffic network is typically characterized by an expensive-to-evaluate objective function. In this paper, we propose two toll level problems (TLPs) integrated with a large-scale simulation-based dynamic traffic assignment (DTA) model of Melbourne, Australia. The first TLP aims to control the pricing zone (PZ) through a time-varying joint distance and delay toll (JDDT) such that the network fundamental diagram (NFD) of the PZ does not enter the congested regime. The second TLP is built upon the first TLP by further considering the minimization of the heterogeneity of congestion distribution in the PZ. To solve the two TLPs, a computationally efficient surrogate-based optimization method, i.e., regressing kriging (RK) with expected improvement (EI) sampling, is applied to approximate the simulation input-output mapping, which can balance well between local exploitation and global exploration. Results show that the two optimal TLP solutions reduce the average travel time in the PZ (entire network) by 29.5% (1.4%) and 21.6% (2.5%), respectively. Reducing the heterogeneity of congestion distribution achieves higher network flows in the PZ and a lower average travel time or a larger total travel time saving in the entire network.

## 1 INTRODUCTION

Congestion pricing as a promising travel demand management strategy has been widely advocated and successfully implemented in cities around the world such as in Singapore, London, Stockholm, and Milan (Gu et al., 2018a). More recently, with various emerging pricing technologies (de Palma and Lindsey, 2011), a few advanced pricing concepts and schemes have emerged including Singapore's electronic road pricing (ERP) system from 2020 onward as a distance-based pricing system, the opt-in distance-based pricing system, OReGo, in Oregon, USA, New York City's Move NY Plan aiming to charge taxis based on both distance and time, and the joint distance- and cordon-based pricing trial in Melbourne, Australia.

Congestion pricing theory has been well established since the seminal studies by Pigou (1920) and Knight (1924). Various first- and second-best pricing models have been proposed as well as solution algorithms (Yang and Huang, 2005). More recently, with the re-theorization of the network fundamental diagram (NFD) or macroscopic fundamental diagram (MFD) (Geroliminis and Daganzo, 2008), a new branch of theory has been enabled that largely facilitates the study, design, and implementation of large-scale pricing. Due to its macroscopic nature, the NFD does not require detailed information about the origin-destination (OD) demand or any individual link in the network, thereby significantly simplifying the modeling and optimization of a large-scale network. See Aboudolas and Geroliminis (2013); Geroliminis et al. (2013); Keyvan-Ekbatani et al. (2012); Ramezani et al. (2015) for NFD-based perimeter control, and Gu et al. (2018c); Simoni et al. (2015); Zheng et al. (2016); Zheng et al. (2012) for NFD-based pricing.

In this paper, we use the NFD to describe the level of congestion in the network and propose a surrogate-based optimization framework to solve two toll level problems (TLPs) in a large-scale network with heterogeneous congestion distribution. A recently developed simulation-based dynamic traffic assignment (DTA) model of Melbourne, Australia is used.

### 1.1 Literature review

The overall toll design problem (TDP) of area-based pricing typically consists of a toll area problem (TAP) (Sumalee, 2004) and a TLP. While few studies tried to solve the TAP and TLP simultaneously (Yang et al., 2002), the majority assumed a predefined pricing zone (PZ) and solved the TLP considering one of the following regimes: (i) zonal (Simoni et al., 2015; Ye et al., 2015), (ii) cordon-based (Liu et al., 2013; Zheng et al., 2016; Zheng et al., 2012), (iii) distance-based (Daganzo and Lehe, 2015; Liu et al., 2017; Meng et al., 2012), and (iv) entry-exit based (Meng and Wang, 2008; Yang et al., 2004b). A joint distance and time toll (JDTT) was recently proposed to overcome the limitation of the distance only toll whereby drivers would intentionally use the shortest paths in the PZ despite being congested (Liu et al., 2014). The defect is that the time toll component as part of the JDTT tends to overcharge drivers as a longer link typically requires a larger travel time despite being uncongested. As such, an improved joint distance and delay toll (JDDT) was proposed (Gu et al., 2018c) and considered in this paper.

The formulation of a TLP in a large-scale dynamic traffic network is often characterized by a computationally expensive objective function, a high-dimensional decision vector, and simulation (if used) noise (Chen et al., 2014). While a stochastic traffic simulator is a source of simulation noise due to different random seed numbers used, many computer experiments involve another type of numerical noise that refers to the random deviations from the expected smooth response (Forrester et al., 2006). Demand uncertainty also contributes to the computational complexity by requiring a significantly larger number of system performance evaluations (Ukkusuri et al., 2007). Therefore, one can hardly devise an analytical method especially without an explicit mathematical model of the system under consideration. Exact gradient methods are no longer applicable, and metaheuristics (Unnikrishnan et al., 2009) are inappropriate given computational concerns. To address this network design problem with an expensive-to-evaluate objective function featuring non-convexity, non-linearity, and non-closed-form, simulation optimization (SO) or simulation-based optimization (SBO) has recently been investigated and advocated (Amaran et al., 2016; Osorio and Bierlaire, 2013). When SO is applied to solve the TLP, existing methods can be classified into two broad categories of feedback control (Gu et al., 2018c; Simoni et al., 2015; Zheng et al., 2016; Zheng et al., 2012) and surrogate-based optimization (Chen et al., 2016; Chen et al., 2014; Chen et al., 2018; Chow and Regan, 2014; Ekström et al., 2016; He et al., 2017).

Feedback control is a classical control strategy aiming to iteratively adjust the control input/decision vector based on the system output/objective function value to meet a desired set point. A pioneering study on feedback pricing control in a large-scale network was recently conducted in an agent-based simulation environment where the NFD was used to describe congestion at the network level (Zheng et al., 2012). The authors applied a discrete integral (I-type) controller to iteratively adjust the cordon toll rate such that the NFD of the PZ does not enter the congested regime. This I-type controller was later improved to a proportional-integral (PI) controller whereby drivers' adaptation to pricing was considered and modeled (Zheng et al., 2016). The authors made a comparison between the two controllers and showed that the latter outperformed the former. Given its superiority, the PI controller was further studied and integrated with a variety of other tolls (Gu et al., 2018c). A similar feedback structure for toll optimization was also proposed by Simoni et al. (2015), although without using a typical feedback controller.

Instead, the authors integrated marginal-cost pricing (MCP) with the NFD to derive their own toll adjustment rule.

Unlike feedback control which resembles a trial-and-error method (Yang et al., 2004a), surrogate-based optimization, also known as response surface method (RSM) or metamodeling, focuses on approximating the simulation input-output mapping using limited function evaluations (Amaran et al., 2016). While surrogate models can be built in local regions to sequentially guide the search for the optimum, global surrogate models from space-filling designs perform better in finding the global optimum (Forrester et al., 2008; Jones et al., 1998). The method, also known as kriging or Gaussian process regression, originates from geostatistics but has become popular in designing and analyzing computer experiments (Sacks et al., 1989). Successful attempts have been made to apply surrogate-based optimization for solving TLPs with different objectives and functional forms of the response surface. Chow and Regan (2014) chose the radial basis function to construct their surrogate model and solved a constrained multi-objective toll optimization problem. However, a comprehensive comparison between different surrogate models revealed that (regressing) kriging with expected improvement (EI) sampling is the best performing surrogate model (Chen et al., 2014; Ekström et al., 2016), which is hence further investigated in a few subsequent studies on toll optimization (Chen et al., 2016; Chen et al., 2018; He et al., 2017).

### 1.2 Objectives and contributions

In this paper, we propose a surrogate-based optimization framework to solve two TLPs in a large-scale dynamic traffic network with heterogeneous congestion distribution. We consider the time-varying JDDT and hence, the entire tolling period is partitioned into several small tolling intervals. Static distance and delay toll rates for each tolling interval are to be optimized, which is different from real-time pricing that produces continuously changing toll rates based on previous measurements.

The objectives of toll pricing can be multiple including but not limited to total travel time minimization (Chen et al., 2014), network speed optimization (Liu et al., 2013), revenue maximization (Saha et al., 2014), network travel time reliability maximization (Chen et al., 2018), and network flow maximization (Zheng et al., 2016). Each of these objectives corresponds to a unique way by which the network is evaluated and hence, different researchers and practitioners may have different preferences. Unlike existing studies on surrogate-based toll optimization, we use the NFD as the network performance indicator based on which two TLPs are formulated. Similar to Simoni et al. (2015); Zheng et al. (2016); Zheng et al. (2012), the first TLP aims to maximize the network flow of the PZ throughout the entire tolling period, but is further subject to a set of toll pattern smoothing constraints, also known as smoothing control constraints (Geroliminis et al., 2013), to help prevent radical changes in the toll rates between adjacent tolling intervals. Therefore, the first TLP is essentially a constrained single-objective optimization problem.

The second TLP is formulated as an extended optimization problem. Given that a heterogeneous distribution of congestion results in a significant hysteresis loop in the NFD (i.e. a key determinant of the shape and scatter of the NFD) causing network unproductivity (Buisson and Ladier, 2009; Knoop and Hoogendoorn, 2013; Mazloumian et al., 2010; Saberi and Mahmassani, 2012; 2013), we introduce another objective into the first TLP to minimize the heterogeneity of congestion distribution in the PZ throughout the entire tolling period. This, to some extent, represents an approach when dealing with large-scale heterogeneous networks (Simoni et al., 2015), as an alternative to network partitioning (Ji and Geroliminis, 2012; Saeedmanesh and Geroliminis, 2016; 2017). Such an objective was previously used to develop a hierarchical perimeter control scheme (Ramezani et al., 2015). Note that while both clustering-based network partitioning and homogeneity control are effective in reducing heterogeneity, heterogeneity itself is an inherent nature of traffic networks that cannot completely disappear (Ramezani et al., 2015). In this paper, instead of solving directly a bi-objective optimization problem, we keep flow maximization as a single objective while reformulating heterogeneity minimization as an additional constraint. A single-objective optimization problem is therefore formulated and to be solved.

The overall contribution of the paper is twofold:

i.  A linear JDDT is investigated and integrated with the NFD to formulate two new high-dimensional TLPs in a large-scale heterogeneously congested dynamic traffic network.
ii. Reducing the heterogeneity of congestion distribution is considered and modeled in the optimization problem, which is shown by the results to help achieve higher network flows.

## 2 PROBLEM FORMULATION

In this section, we first formulate the JDDT and the resultant generalized travel cost function to be used in the simulation model for vehicle routing and path assignment. We then develop two TLPs to optimize the toll rates such that the network performance objective(s) is achieved.

### 2.1 Joint distance and delay toll (JDDT)

Consider a network $G = (N, A)$ where $N$ is the set of nodes and $A$ is the set of directed links. With a predefined pricing cordon, network $G$ is partitioned into a PZ $G_p = (N_p, A_p)$ and a peripheral sub-network $G_{np} = (N_{np}, A_{np})$. Table 1 summarizes the notation used to formulate the JDDT.

Let $\boldsymbol{\tau} = [v_1, \ldots, v_m, \omega_1, \ldots, \omega_m]^T$ be the toll decision vector for the $m$ tolling intervals where $v_h$ and $\omega_h$ are the distance toll rate and the delay toll rate for the $h$-th tolling interval, respectively. We assume that the distance and the delay toll functions are both linear with respect to the distance traveled and the delay experienced within the PZ, respectively. Here, delay is defined as the difference between the actual simulated travel time and the free-flow travel time. This assumption enables the link additive property resulting in a link-based approach to pricing system modeling and optimization. If non-linearity is assumed, a path-based approach should be pursued (Liu et al., 2017). Let $l_r^{od}(h)$ and $t_r^{od}(h)$ be the distance traveled and the time spent within the PZ for path $r \in R^{od}$ during the $h$-th tolling interval, respectively:

**Table 1**
Notation used to formulate the JDDT

| Notation | Interpretation |
|---|---|
| $W$ | The set of OD pairs where $O \subset N$ is the set of origins and $D \subset N$ is the set of destinations, i.e. $W = \{(o,d)\|o \in O, d \in D\}$ |
| $R^{od}$ | The set of paths between an OD pair $(o,d) \in W$ |
| $m$ | Total number of tolling intervals |
| $l_a$ | Length of link $a \in A$ |
| $t_a(h)$ | Average travel time on link $a \in A$ during the $h$-th tolling interval |
| $\delta_{a,r}^{od}$ | $\delta_{a,r}^{od} = 1$ if path $r \in R^{od}$ includes link $a$, otherwise $\delta_{a,r}^{od} = 0$ |

$$l_r^{od}(h) = \sum_{a \in A_p} l_a \delta_{a,r}^{od} \tag{1}$$

$$t_r^{od}(h) = \sum_{a \in A_p} t_a(h) \delta_{a,r}^{od} \tag{2}$$

where $r \in R^{od}, (o,d) \in W, h \in (1,2,\ldots,m)$. The distance toll component and the delay toll component, $\varphi_r^{od}(h)$ and $\phi_r^{od}(h)$, for path $r \in R^{od}$ during the $h$-th tolling interval are therefore defined and calculated, respectively, as

$$\varphi_r^{od}(h) = v_h \sum_{a \in A_p} l_a \delta_{a,r}^{od} \tag{3}$$

$$\phi_r^{od}(h) = \omega_h \sum_{a \in A_p} \left(t_a(h) - t_a^f\right) \delta_{a,r}^{od} \tag{4}$$

where $t_a^f$ is the free-flow travel time on link $a \in A_p$. The generalized travel cost function, $V_r^{od}(h)$, for path $r \in R^{od}$ during the $h$-th tolling interval is expressed as:

$$V_r^{od}(h) = \sum_{a \in A} t_a(h) \delta_{a,r}^{od} + \frac{\varphi_r^{od}(h) + \phi_r^{od}(h)}{VTT} \tag{5}$$

where VTT is drivers' average value of travel time. In the simulation model, Equation (5) is integrated with the C-logit stochastic route choice model for path assignment (Cascetta et al., 1996).

## 2.2 Single-objective toll level problem (TLP)

In the single-objective TLP, we aim to optimize the time-varying JDDT such that the NFD of the PZ does not enter the congested regime. Assuming the critical network density of the PZ does not change significantly before and after pricing, the single-objective TLP is formulated as follows:

$$\min_{\boldsymbol{\tau}\in\Omega} \mathrm{E}\left[\frac{1}{m}\sum_{h=1}^{m}|\bar{K}_h - K_{\mathrm{cr}}|\right] \quad (6)$$

s.t.

$$|v_h - v_{h+1}| \leq \alpha, \quad h = 1,2,\dots,m \quad (7)$$
$$|\omega_h - \omega_{h+1}| \leq \beta, \quad h = 1,2,\dots,m \quad (8)$$
$$\bar{K}_h = DTA(\boldsymbol{\tau}), \quad h = 1,2,\dots,m \quad (9)$$
$$\Omega = \{\boldsymbol{\tau}|\boldsymbol{\tau}_{\min} \leq \boldsymbol{\tau} \leq \boldsymbol{\tau}_{\max}\} \quad (10)$$

where $\bar{K}_h$ is the average network density of the PZ during the $h$-th tolling interval, $K_{\mathrm{cr}}$ is the critical network density of the PZ identified from the non-tolling NFD, $\alpha$ and $\beta$ are the toll pattern smoothing parameters for the distance and the delay toll rates, respectively, and $\Omega$ is the feasible set of toll rates with $\boldsymbol{\tau}_{\min}$ being the lower bound and $\boldsymbol{\tau}_{\max}$ being the upper bound.

Given a stochastic traffic simulator, the objective function aims to minimize the expected average of the absolute difference between $\bar{K}_h$ and $K_{\mathrm{cr}}$ for the $m$ tolling intervals. As such, the PZ is driven to achieve the maximum network flow throughout the entire tolling period. To calculate the expectation, one can readily apply fixed-number sample path optimization, also known as sample average approximation (Amaran et al., 2016). However, in the presence of a computationally expensive objective function, the sample size is usually restricted to a small number just to reduce the effect of noise rather than pursuing a complete noise filter. To handle simulation noise more efficiently, one can apply variable-number sample path optimization (He et al., 2017). Equations (7) and (8) are the toll pattern smoothing constraints used to ensure that the optimal toll rates do not fluctuate unduly between adjacent tolling intervals, and that we obtain a smooth optimal toll pattern. It would be practically infeasible to introduce a radically changing pricing scheme considering drivers' adaptivity and system stability. Equation (9) represents the black-box function of the simulation-based DTA model to which $\boldsymbol{\tau}$ is input to obtain $\bar{K}_h$.

## 2.3 Bi-objective toll level problem (TLP)

In the bi-objective TLP, we further aim to minimize the heterogeneity of congestion distribution in the PZ throughout the entire tolling period. To this end, we introduce the spatial spread of density, $\gamma$, as the square root of the weighted variance of link densities (Knoop and Hoogendoorn, 2013) to quantify the heterogeneity of congestion distribution:

$$\gamma = \sqrt{\frac{\sum_i l_i n_i (k_i - K)^2}{\sum_i l_i n_i}}, \qquad K = \frac{\sum_i k_i l_i n_i}{\sum_i l_i n_i} \quad (11)$$

where $k_i$ is the average density of link $i$ over the observation period, and $l_i$ and $n_i$ are the length and the number of lanes of link $i$, respectively. When $l_i n_i$ is the same for every link $i$ in the network, Equation (11) reduces to the standard deviation measure (Mahmassani et al., 2013; Mazloumian et al., 2010; Saberi and Mahmassani, 2012). The spatial spread of density naturally increases with a growing accumulation – an increase in vehicles entering the PZ inevitably generates a higher spatial spread of density later in time as these vehicles continue their trips. Therefore, a better quantification of the heterogeneity of congestion distribution according to Simoni et al. (2015) are the positive deviations from the natural increment represented by the lower envelope of the spread-accumulation relationship, which is termed the deviation from spread, $\Delta$:

$$\Delta = \gamma - \gamma(K) \quad (12)$$

where $\gamma(K) = aK^3 + bK^2 + cK$ is an assumed third-order polynomial function fitted to the lower envelope of the spread-accumulation relationship with coefficients $a$, $b$, and $c$ to be estimated. Note that the lower envelope of the spread-accumulation relationship corresponds to the upper envelope of the NFD, because the least possible heterogeneity of congestion distribution contributes to the highest possible network flow for a certain network density. Also note that the fitted function only serves as a mathematical approximation and hence does not necessarily represent the best functional form.

The bi-objective TLP is therefore formulated as follows:

$$\min_{\boldsymbol{\tau}\in\Omega} \mathrm{E}\left[\frac{1}{m}\sum_{h=1}^{m}|\bar{K}_h - K_{\mathrm{cr}}|\right] \quad (13)$$

$$\min_{\boldsymbol{\tau}\in\Omega} \mathrm{E}\left[\frac{1}{m}\sum_{h=1}^{m}\bar{\Delta}_h\right] \quad (14)$$

s.t.

$$|v_h - v_{h+1}| \leq \alpha, \quad h = 1,2,\dots,m \quad (15)$$
$$|\omega_h - \omega_{h+1}| \leq \beta, \quad h = 1,2,\dots,m \quad (16)$$
$$\bar{K}_h = \mathrm{DTA}(\boldsymbol{\tau}), \quad h = 1,2,\dots,m \quad (17)$$
$$\Omega = \{\boldsymbol{\tau}|\boldsymbol{\tau}_{\min} \leq \boldsymbol{\tau} \leq \boldsymbol{\tau}_{\max}\} \quad (18)$$

where $\bar{\Delta}_h$ is the average deviation from spread of the PZ during the $h$-th tolling interval. The unique feature of Problem (13-18) is that we know a priori that both objective functions have a lower bound of zero (because Equation (13) involves an absolute value operator and we use the lower envelope to calculate $\Delta$ in Equation (12)), although being too ideal to achieve. While we could solve the bi-objective optimization problem as it is, we could alternatively utilize this unique feature by keeping Equation (13) as the single objective and reformulating Equation (14) as an additional constraint. The original bi-objective TLP is therefore transformed into the following single-objective optimization problem:

$$\min_{\boldsymbol{\tau} \in \Omega} \mathrm{E}\left[\frac{1}{m}\sum_{h=1}^{m}|\bar{K}_h - K_{\mathrm{cr}}|\right] \quad (19)$$

s.t.

$$\mathrm{E}\left[\frac{1}{m}\sum_{h=1}^{m}\bar{\Delta}_h\right] \leq \Delta_{\max} \quad (20)$$

$$|v_h - v_{h+1}| \leq \alpha, \quad h = 1,2,\ldots,m \quad (21)$$
$$|\omega_h - \omega_{h+1}| \leq \beta, \quad h = 1,2,\ldots,m \quad (22)$$
$$\bar{K}_h = DTA(\boldsymbol{\tau}), \quad h = 1,2,\ldots,m \quad (23)$$
$$\Omega = \{\boldsymbol{\tau}|\boldsymbol{\tau}_{\min} \leq \boldsymbol{\tau} \leq \boldsymbol{\tau}_{\max}\} \quad (24)$$

where $\Delta_{\max}$ is a constraint limit to ensure that the heterogeneity of congestion distribution is below a certain threshold. How to determine $\Delta_{\max}$ depends on the network performance under consideration, and hence is case-specific. In general, we apply trail-and-error by utilizing the knowledge of the network performance under the non-tolling scenario as well as the results of the single-objective TLP – $\Delta_{\max}$ should be smaller than $\mathrm{E}\left[\frac{1}{m}\sum_{h=1}^{m}\bar{\Delta}_h\right]$ from solving the single-objective TLP and larger than that under the non-tolling scenario if the NFD enters the congested regime without experiencing network recovery. While we could alternatively keep Equation (14) as the objective and reformulate Equation (13) as the constraint, we choose not to so as to be consistent with the formulation of Problem (6-10).

Note that one could also construct a weighted average of the two objectives to formulate a single-objective optimization problem (Chen et al., 2016), but the key question is how to choose the relative weight (perhaps through a sensitivity analysis). There is no need to re-organize the complex constraints as part of the objective function, e.g., through a penalty method or Lagrangian relaxation. This is because the surrogate-based optimization method to be introduced and applied can well handle complex constraints.

## 3 SURROGATE-BASED OPTIMIZATION

To solve the two TLPs, a surrogate-based optimization method is applied featuring regressing kriging (RK) with expected improvement (EI) sampling (see Figure 1). To construct the starting surrogate model, a few initial sample points need to be generated through a space-filling design of experiments (DOE), for each of which a network simulation is performed to evaluate the objective function. The constructed surrogate model is further subject to adding infill sample points via EI sampling until model validation is passed. Here, accuracy measures how well RK predicts while convergence indicates whether there is still room for augmenting RK to improve the objective function value. It is possible that an infill sample point does not improve the current best solution since the constructed response surface used for calculating and maximizing the EI is an approximation to the true unknown response surface. As a result, the surrogate-based optimum approximates the true unknown optimum. For most practical applications with tight computational budgets, a maximum number of iterations or function evaluations is usually reached first before a good convergence is achieved suggesting that the final solution is not assumed to be optimal (Amaran et al., 2016; Osorio and Punzo, 2019).

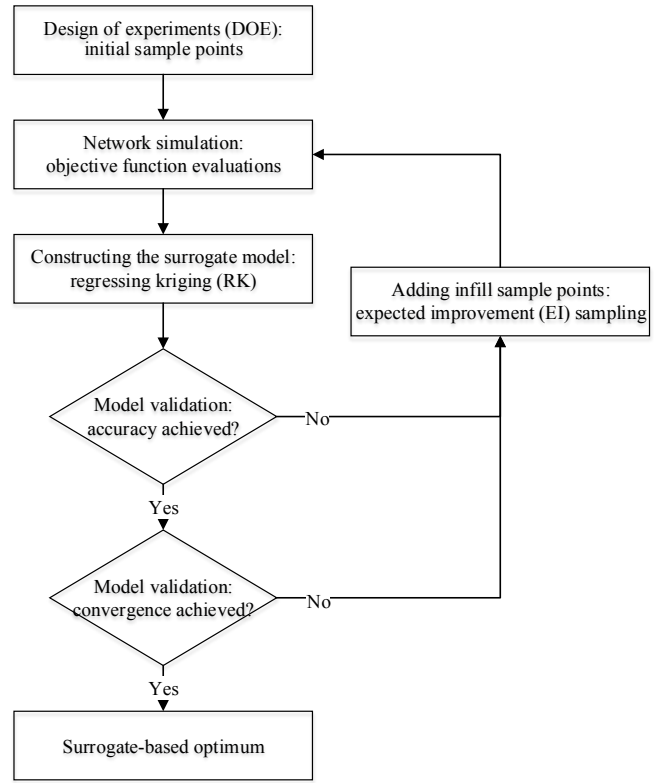

**Figure 1** Flowchart representation the surrogate-based optimization method

### 3.1 Design of Experiments (DOE)

Since DOE aims to provide an initial set of sample points to construct the starting surrogate model, the space-filling property is desirable as the resultant sample points are spread as uniformly as possible over the entire feasible domain. Latin Hypercube Sampling (LHS) is a space-filling DOE whereby each problem dimension is stratified into an equal number of intervals from which points are uniformly sampled. As such, there is no overlap in LHS when mapping the multi-dimensional sample points into each dimension. To achieve the maximum uniformity or space-fillingness of an LHS plan, one can apply maximin LHS to maximize the minimum distance between all the sample points by generating and evaluating a set of candidate plans (Forrester et al., 2008). According to Ekström et al. (2016), at least $2m + 1$ sample points are required to construct the starting surrogate model where $2m$ is the problem dimension, i.e., the size of the toll decision vector. Due to the much higher dimension of our toll optimization problem, we choose the size of the initial set of sample points to be $2(2m + 1)$. A few additional sample points can also be considered as part of the initial plan such as the corner and center points of the design space.

### 3.2 Regressing kriging (RK)

Kriging metamodeling is a stochastic process approach originating from Bayesian reasoning where the output of a deterministic computer experiment is modeled as a realization of a stochastic process. The well-known ordinary kriging model assumes an unknown constant mean $\mu$ of the response surface $y(\boldsymbol{\tau})$ where $\boldsymbol{\tau}$ is the toll decision vector, and a zero-mean second-order stationary Gaussian process $Z$, i.e., $y(\boldsymbol{\tau}) = \mu + Z(\boldsymbol{\tau})$ where $E[Z(\boldsymbol{\tau})] = 0$. Here, $y(\boldsymbol{\tau})$ equates to $E\left[\frac{1}{m}\sum_{h=1}^{m}|\bar{K}_h - K_{cr}|\right]$ in Equation (19). The covariance function of $Z$ between any two points $\boldsymbol{\tau}^{(i)}$ and $\boldsymbol{\tau}^{(j)}$ is $\sigma^2 \psi(\boldsymbol{\tau}^{(i)}, \boldsymbol{\tau}^{(j)})$ where $\sigma^2$ is the process variance and $\psi(\cdot)$ is the Gaussian correlation function depending on the distance between $\boldsymbol{\tau}^{(i)}$ and $\boldsymbol{\tau}^{(j)}$ only, i.e., $\psi(\boldsymbol{\tau}^{(i)}, \boldsymbol{\tau}^{(j)}) = \exp\left(-\sum_{l=1}^{2m} \theta_l \left(\boldsymbol{\tau}_l^{(i)} - \boldsymbol{\tau}_l^{(j)}\right)^2\right)$ where $\boldsymbol{\theta}$ is a vector of scaling coefficients that allows for varying impacts of different dimensions on the correlation function and $2m$ is the dimension of the toll optimization problem formulated in Section 2. The correlation matrix $\boldsymbol{\Psi}$ is constructed with the $(i,j)$-th element being $\psi(\boldsymbol{\tau}^{(i)}, \boldsymbol{\tau}^{(j)})$.

The ordinary kriging model is an interpolation method that constructs the response surface by passing through all the sample points. When computer simulations display the numerical noise, i.e., the output tend to have a random scatter about a smooth trend rather than lying on it, the interpolating kriging model may exhibit overfitting without being able to tolerate data fluctuations (Forrester et al., 2006). The solution is to allow the kriging model not to interpolate but to regress the sample points, which is achieved by adding a regularization constant, $\lambda$, to the diagonal of the correlation matrix. That is, $\mathbf{R} = \boldsymbol{\Psi} + \lambda \mathbf{I}$ where $\mathbf{R}$ is known as the regressing correlation matrix and $\mathbf{I}$ is an identity matrix of the same dimension. The resultant model is commonly known as RK (Chen et al., 2014; Forrester et al., 2006; He et al., 2017).

Given the assumption of a Gaussian process, the kriging predictor and the prediction error can be obtained by maximizing the augmented log-likelihood as shown in Forrester et al. (2006):

$$\hat{y}(\boldsymbol{\tau}^*) = \hat{\mu} + \boldsymbol{\psi}^T \mathbf{R}^{-1}(\mathbf{y} - \mathbf{1}\hat{\mu}) \qquad (25)$$
$$\hat{s}^2(\boldsymbol{\tau}^*) = \hat{\sigma}^2\left(1 + \hat{\lambda} - \boldsymbol{\psi}^T \mathbf{R}^{-1} \boldsymbol{\psi}\right) \qquad (26)$$

where $\boldsymbol{\psi} = \left[\psi(\boldsymbol{\tau}^*, \boldsymbol{\tau}^{(1)}), \psi(\boldsymbol{\tau}^*, \boldsymbol{\tau}^{(2)}), \ldots, \psi(\boldsymbol{\tau}^*, \boldsymbol{\tau}^{(n)})\right]^T$ is the correlation vector between the new point $\boldsymbol{\tau}^*$ and all the existing sample points. Equations (25) and (26) clearly show the essence of RK or Bayesian reasoning – the prediction is modeled as a distribution with a mean and a variance (namely a random variable), rather than being deterministic. In general, we are more confident in the prediction if the variance is low, and vice versa. To obtain the maximum likelihood estimates (MLEs) of the parameters involved in the log-likelihood function, i.e., $\hat{\mu}$, $\hat{\sigma}^2$, $\hat{\lambda}$, and $\hat{\boldsymbol{\theta}}$, we apply the genetic algorithm (GA).

### 3.3 EI sampling

When kriging is used to approximate the simulation input-output mapping, additional infill sample points are required to enhance the constructed response surface. In general, there are two categories of infill strategies (Ekström et al., 2016):

i. One-stage infill strategies which search for infill sample points according to a certain merit function, e.g., maximizing the minimum distance between all the sample points, without using information about the constructed response surface
ii. Two-stage infill strategies which search for infill sample points by utilizing the constructed response surface

We choose a two-stage infill strategy given its self-learning mechanism – the new response surface is iteratively augmented based on its predecessor. Specifically, we apply a global optimal infill strategy known as EI sampling as opposed to a suboptimal infill strategy that balances poorly between exploring unvisited regions and exploiting visited regions (Chen et al., 2014). While trying to locate infill sample points that lead to low predictor values for a minimization problem, EI sampling also considers uncertainty about the constructed response surface as reflected by the prediction error. In regions with few sample points, although the current prediction may not be promising, the error is likely to be high suggesting a good opportunity to improve the current best solution by adding infill sample points. Therefore, as a global search method, EI sampling can balance well between local exploitation and global exploration (Forrester et al., 2008). To solve the EI maximization problem, we use the GA.

#### 3.3.1 Unconstrained EI sampling

Unconstrained EI sampling only considers maximizing the EI of the objective when adding infill sample points. At each iteration, a new response surface approximating the objective function in Equation (6) is constructed using all the exiting sample points. This enables us to predict for any point in the design space using Equations (25) and (26). Our goal is then to find the point that maximizes the EI because a lager EI equates to a greater probability of improving the current best solution, and to evaluate this identified point through the simulation model for further augmenting the response surface. Let $y_{min}$ denote the best observed objective function value so far. The improvement at a new infill sample point $\boldsymbol{\tau}^*$ is defined as $I(\boldsymbol{\tau}^*) = \max(y_{min} - y(\boldsymbol{\tau}^*), 0)$. The EI at this point hence reads $E[I(\boldsymbol{\tau}^*)] = E[\max(y_{min} - y(\boldsymbol{\tau}^*), 0)]$. When $\hat{s}^2(\boldsymbol{\tau}^*) = 0$, $E[I(\boldsymbol{\tau}^*)] = 0$; when $\hat{s}^2(\boldsymbol{\tau}^*) > 0$, given the assumption of a Gaussian process,

$$E[I(\boldsymbol{\tau}^*)] = \frac{1}{\sqrt{2\pi \hat{s}^2(\boldsymbol{\tau}^*)}} \int_{-\infty}^{y_{min}} (y_{min} - u) \exp\left(-\frac{(u - \hat{y}(\boldsymbol{\tau}^*))^2}{2\hat{s}^2(\boldsymbol{\tau}^*)}\right) du \quad (27)$$

When using ordinary kriging, both the prediction error and the EI stay at zero for all the existing sample points. It is therefore impossible to add an infill sample point that has already been sampled. However, when using RK, $\hat{s}^2(\boldsymbol{\tau}^*) = 0$ does not hold at an existing sample point resulting in the possibility of maximizing $E[I(\boldsymbol{\tau}^*)]$ at a previously sampled point. To prevent RK from getting trapped at an existing sample point, Forrester et al. (2006) proposed a reinterpolation MLE of $\sigma^2$:

$$\hat{\sigma}_{ri}^2 = \frac{(\mathbf{y} - \mathbf{1}\hat{\mu})^T \mathbf{R}^{-1} \boldsymbol{\Psi} \mathbf{R}^{-1} (\mathbf{y} - \mathbf{1}\hat{\mu})}{n} \quad (28)$$

The reinterpolation prediction error reads $\hat{s}_{ri}^2(\boldsymbol{\tau}^*) = \hat{\sigma}_{ri}^2(1 - \boldsymbol{\psi}^T \mathbf{R}^{-1} \boldsymbol{\psi})$. Now, $\hat{s}_{ri}^2(\boldsymbol{\tau}^*) = 0$ holds for all the existing sample points for which $E[I_{ri}(\boldsymbol{\tau}^*)] = 0$. When $\hat{s}_{ri}^2(\boldsymbol{\tau}^*) > 0$ and assuming $y(\boldsymbol{\tau}^*) \sim N(\hat{y}(\boldsymbol{\tau}^*), \hat{s}_{ri}^2(\boldsymbol{\tau}^*))$, $E[I_{ri}(\boldsymbol{\tau}^*)]$ can be calculated through Equation (27) by simply replacing $\hat{s}^2(\boldsymbol{\tau}^*)$ with $\hat{s}_{ri}^2(\boldsymbol{\tau}^*)$.

3.3.2 Constrained EI sampling

While maximizing the EI of the objective, constrained EI sampling further considers the impact of the constraint on adding infill sample points. At each iteration, apart from constructing a new response surface to approximate the objective function in Equation (19) using all the exiting sample points, we also construct a new response surface to approximate the constraint in Equation (20). By calculating the EI of the objective and the probability of not violating the constraint, and maximizing their product, we are able to find the point that potentially improves the current best solution most while satisfying the constraint.

Let $c(\boldsymbol{\tau})$ denote the response surface of the constraint to be no greater than $\Delta_{max}$, i.e., $c(\boldsymbol{\tau})$ equates to $E\left[\frac{1}{m}\sum_{h=1}^{m} \bar{\Delta}_h\right]$ in Equation (20). The constrained improvement at a new infill sample point is defined as

$$CI(\boldsymbol{\tau}^*) = \begin{cases} I(\boldsymbol{\tau}^*), & c(\boldsymbol{\tau}^*) \leq \Delta_{max} \\ 0, & c(\boldsymbol{\tau}^*) > \Delta_{max} \end{cases} \quad (29)$$

If the constraint is violated at $\boldsymbol{\tau}^*$, i.e., $c(\boldsymbol{\tau}^*) > \Delta_{max}$, $CI(\boldsymbol{\tau}^*)$ is zero even if $y_{min} - y(\boldsymbol{\tau}^*)$ is large. The constrained EI hence reads $E[CI(\boldsymbol{\tau}^*)] = E[I(\boldsymbol{\tau}^*)]P[c(\boldsymbol{\tau}^*) \leq \Delta_{max}]$ where $P[c(\boldsymbol{\tau}^*) \leq \Delta_{max}]$ is the probability of not violating the constraint. The constrained EI is large only if the EI of the objective and the probability of not violating the constraint are both large. With reinterpolation, we end up with $E[CI_{ri}(\boldsymbol{\tau}^*)] = E[I_{ri}(\boldsymbol{\tau}^*)]P_{ri}[c(\boldsymbol{\tau}^*) \leq \Delta_{max}]$ where

$$P_{ri}[c(\boldsymbol{\tau}^*) \leq \Delta_{max}]$$
$$= \frac{1}{\sqrt{2\pi \hat{s}_{cri}^2(\boldsymbol{\tau}^*)}} \int_{-\infty}^{\Delta_{max}} \exp\left(-\frac{(u - \hat{c}(\boldsymbol{\tau}^*))^2}{2\hat{s}_{cri}^2(\boldsymbol{\tau}^*)}\right) du \quad (30)$$

### 3.4 Model validation

To validate the accuracy of the constructed surrogate model, one option is to select a few additional sample points to form a test set based on which the observed and predicted objective function values are compared. The training set obviously includes the initial sample points and those added as infill sample points. This option, however, is not desirable particularly when concern about the extra computational effort prevails. A better option which has been adopted in a few relevant studies (Chen et al., 2014; Ekström et al., 2016) is to leave out one observation and predict it based on the remaining observations. This procedure is commonly known as the leave-one-out cross validation (CV) which requires no additional sample points to validate the accuracy of the model.

With a total of $n$ observations, the leave-one-out CV is repeated $n$ times and each time it produces a cross-validated prediction $\hat{y}_{-i}(\boldsymbol{\tau}^{(i)})$ for the corresponding observation $y(\boldsymbol{\tau}^{(i)})$. While common measures of effectiveness (MOEs) can be calculated to reflect the prediction accuracy, they are inappropriate for evaluating the surrogate model as the prediction at any point is a normally distributed random variable rather than a scalar (Chen et al., 2014). Knowing that, along with the cross-validated prediction, we also obtain a cross-validated standard error $\hat{s}_{-i}(\boldsymbol{\tau}^{(i)})$, we can calculate the 99.7% confidence interval for each $y(\boldsymbol{\tau}^{(i)})$ using the prediction plus or minus three standard errors (Jones et al., 1998). Alternatively, we can calculate $\frac{y(\boldsymbol{\tau}^{(i)}) - \hat{y}_{-i}(\boldsymbol{\tau}^{(i)})}{\hat{s}_{-i}(\boldsymbol{\tau}^{(i)})}$ to obtain a standardized cross-validated residual, the value of which should be lying roughly within $[-3,3]$ for an accurate surrogate model. Unfortunately, there is no clean proof of convergence for the surrogate-based optimization method. A practical technique is to track the convergence history of the EI as we will show in Section 4.

### 4 RESULTS AND DISCUSSION

#### 4.1 Experimental setup

In this paper, we employ a recently developed large-scale simulation-based DTA model of Melbourne, Australia

(Shafiei et al., 2018). The model is deployed in AIMSUN with time-varying commuting demand during the 6-10 AM peak period. By using multi-source big traffic data, major supply (Gu et al., 2016; 2018b) and demand (Shafiei et al., 2018) input to the model has been calibrated and validated. Figure 2(a) shows the extracted sub-network from the greater Melbourne area model. There are 4,375 links and 1,977 nodes in the entire network, and 282 links and 91 nodes in the PZ represented by the inner rectangle. During the simulation, path assignment is calculated every 15 minutes using the C-logit stochastic route choice model assuming VTT = $15/h. 30% of drivers are assumed adaptive having access to real-time information and hence, can update their shortest paths enroute at the beginning of every path assignment interval using information of traffic conditions from the previous interval. Note that the simulation model does not explicitly account for how information is communicated between vehicles and the network infrastructure. A sensitivity analysis on the percentage of adaptive driving is provided in Appendix B.

To determine the critical network density and the tolling period, we run simulation without pricing and show the density time series and the simulated NFDs of the PZ in Figure 2(b) and (c), respectively. Results suggest that we set $K_{cr}$ at 25 vpkmpl which leads to a 2-h tolling period between 8 (a few minutes before the network density reaches the 25 threshold and the network becomes congested) and 10 AM.

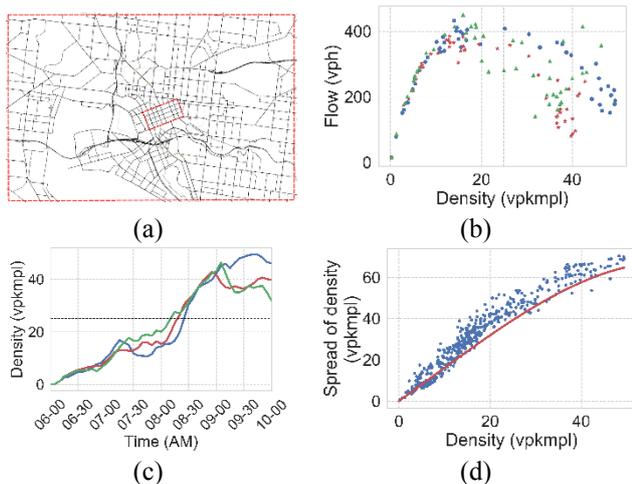

**Figure 2** (a) Extracted sub-network from the greater Melbourne area model and the simulation results of the PZ under the non-tolling scenario: (b) density time series, (c) simulated NFDs, and (d) spread-accumulation relationship

To demonstrate the capability of the surrogate-based optimization framework in dealing with high-dimensional problems, we use a 15-min duration and partition the entire tolling period into 8 small tolling intervals. Hence a total of 16 toll decision variables are to be optimized, 8 of which are distance toll rates and the other 8 are delay toll rates.

Accordingly, in the maximin LHS plan, the total number of the initial sample points is 37. When applying the surrogate-based optimization framework, we allow a maximum of 100 iterations, i.e., the total number of sample points is 100 with 63 infill sample points. Without loss of generality, $\tau_{\min}$ and $\tau_{\max}$ are set at $[0,...,0,0,...,0]^T$ and $[1,...,1,15,...,15]^T$, respectively, and $\alpha$ and $\beta$ are set at $\frac{1}{3}(1-0) \approx 0.33$ and $\frac{1}{3}(15-0) = 5$, respectively. A sensitivity analysis on $\alpha$ and $\beta$ is provided in Appendix A. To estimate $a$, $b$, and $c$, we run ten replications without pricing and obtain the following fitted functional form, $\gamma(K) = -0.0002032K^3 + 0.004432K^2 + 1.587K$, which is also shown in Figure 2(d). Although being estimated under the non-tolling scenario, $\gamma(K)$ is applicable to different pricing scenarios because it captures the highest possible network flow for a certain network density, which can be considered as the invariant "capacity" of the network corresponding to this network density.

### 4.2 Solving the single-objective toll level problem (TLP)

Figure 3(a) validates the accuracy of the constructed surrogate model with 100 sample points. The model accuracy is sufficiently achieved with 98 standardized cross-validated residuals lying within $[-3,3]$. One outlier corresponds to the non-tolling scenario with $\tau = \tau_{\min} = [0,0,...,0,0,0,...,0]^T$. Since the non-tolling network produces the highest objective function value, the surrogate model makes little effort exploring the region surrounding the non-tolling sample point where the prediction becomes poor, as expected. Figure 3(b) illustrates the convergence history of the EI. Although, due to the heuristic nature of the method, intermittent peaks representing possible significant improvements in the objective function value are observed, the overall trend of the change as represented by the average curve (averaged every four consecutive points) displays a relatively smooth convergence pattern towards zero. This implies that, at the end of optimization, the surrogate model is unable to locate a new solution that significantly improves the current best solution and hence, we can terminate the algorithm with confidence.

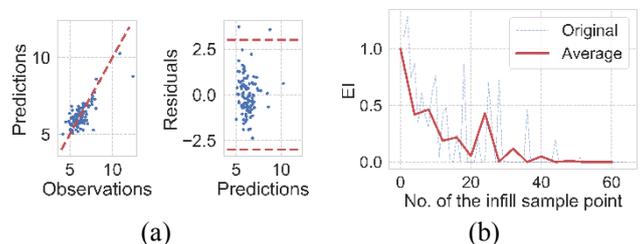

**Figure 3** Solving the single-objective TLP: (a) validating the accuracy of the constructed surrogate model, and (b) convergence history of the EI

The solution to the single-objective TLP is shown in Figure 4(a). The changes in the distance and delay toll rates between adjacent tolling intervals are clearly bounded by the

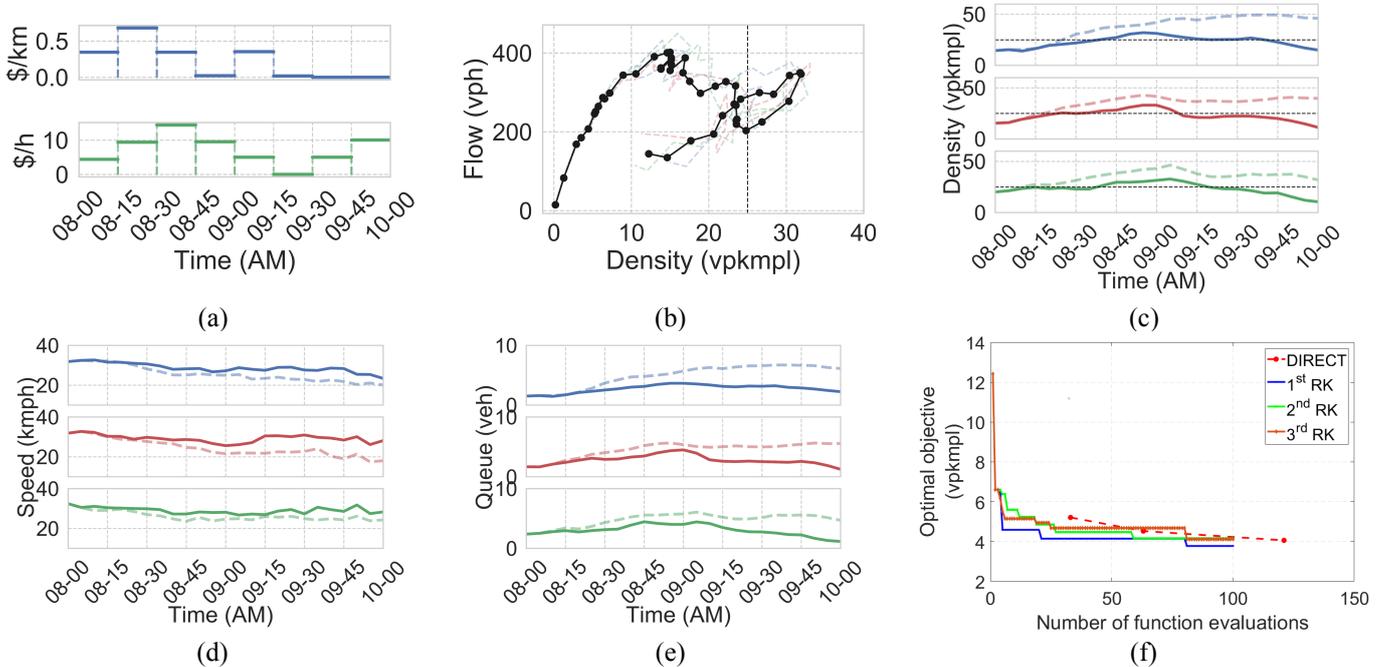

**Figure 4** (a) Optimal distance and delay toll patterns from solving the single-objective TLP, (b) simulated NFDs for all replications (dashed lines) and the averaged NFD (solid line) of the PZ, (c-e) time series of density, speed, and queue of the PZ(dashed/solid lines represent the before-/after-pricing scenario), and (f) changes in the optimal objective function value as the number of function evaluations increases for both RK and DIRECT toll pattern smoothing constraints in Equations (7) and (8), respectively. Figure 4(b) shows the simulated averaged NFD of the PZ after applying the optimal toll rates. As expected, the congested regime of the NFD that appears and remains until the end of simulation under the non-tolling scenario no longer exists and is substituted by a combination of a (near-)capacity regime and a clock-wise hysteresis loop.

When the network is unloading, namely recovering from congestion, the distribution of congestion tends to be more heterogeneous than when the network is loading, because the congested area often clear slowly and are fragmented during network unloading/recovery. The heterogeneous distribution of congestion inevitably reduces the network flow giving rise to a clockwise hysteresis loop in the NFD (Gayah and Daganzo, 2011; Geroliminis and Sun, 2011; Saberi and Mahmassani, 2013). An interesting observation out of the comparison is that, compared with the non-tolling NFD, the tolling NFD undergoes a network flow drop immediately after the implementation of pricing, which, in part, contributes to the hysteresis loop in the NFD. This drop, as was also observed in Gu et al. (2018c), results from the reduced inflow or demand to the PZ due to the presence of pricing. An extreme and apparently unrealistic scenario is that we implement an exceptionally high toll price whereby no one would enter the PZ. With such demand dropping sharply to zero, the hysteresis loop in the NFD is amplified most significantly (Mahmassani et al., 2013). A complete elimination of the network flow drop is too ideal and perhaps only possible with an extremely smooth toll pattern starting from zero, i.e. a very slow-varying toll. Figure 4(c-e) show, respectively, the density, speed, and queue time series of the PZ under the optimal tolling scenario in comparison with those under the non-tolling scenario. It is evident and consistent across different replications that pricing has brought significant performance improvement to the PZ represented by the area in between the two curves.

The computational efficiency of the surrogate-based optimization method has been highlighted, e.g., in Chen et al. (2014); Chow et al. (2010) in comparison with the GA as a random search optimization method. To further show the computational efficiency of RK, we apply another global optimization method termed DIviding RECTangles (DIRECT) to solve the single-objective TLP in comparison with RK. DIRECT is a deterministic method originating from Lipschitzian optimization. It works by iteratively partitioning the search space into multiple hyperrectangles and identifying what are called the potentially optimal hyperrectangles for further partitioning. Details of DIRECT and its mathematical properties can be found in Jones et al. (1993). However, DIRECT cannot be applied directly due to the presence of the toll pattern smoothing constraints in the TLP. We therefore integrate DIRECT with the penalty function method (Bazaraa et al., 2013) to transform the original constrained optimization problem into an unconstrained one. To be consistent with RK, we terminate DIRECT when the number of function evaluations exceeds 100. We perform RK for a total of three runs given its stochastic nature arising from LHS and the GA, but only one run of DIRECT is performed because of its deterministic nature. With 121 function evaluations, the optimal objective function value of DIRECT is 4.06 vpkmpl, while with 100 function evaluations, the optimal objective function values of RK are 3.78, 4.17, and 4.11 vpkmpl, respectively. Although the solution quality of RK and DIRECT is similar, the latter requires a larger number of function evaluations as shown in Figure 4(f). Since DIRECT identifies all the potentially optimal hyperrectangles and evaluates their center points at each iteration, there are only three points along the curve representing the three iterations before termination. The three iterations perform 33, 30, and 58 function evaluations, respectively, which add up to the 121 function evaluations. Within 50 function evaluations, RK manages to find a solution that is very close to the final optimum, whereas DIRECT needs 63 function evaluations to achieve roughly the same level of optimality. RK manages to locate the final optimum with about 80 function evaluations, whereas DIRECT needs a much larger 121 function evaluations. Since each function evaluation (i.e., one simulation run) takes on average 15-20 minutes, RK can save more than 10 hours of simulation time compared with DIRECT and hence, has better performance in terms of the computational efficiency.

### 4.3 Solving the bi-objective toll level problem (TLP)

When solving the bi-objective TLP, we set $\delta_{max}$ at 8 vpkmpl in Equation (20). Figure 5(a) validates the accuracy of the constructed surrogate model with 100 sample points. As with Figure 3(a), there are 98 well-predicted sample points plus two outliers. One of the outliers still corresponds to the non-tolling scenario with $\boldsymbol{\tau} = \boldsymbol{\tau}_{min} = [0, ..., 0, 0, ..., 0]^T$, while the other outlier corresponds to the "full" tolling scenario with $\boldsymbol{\tau} = \boldsymbol{\tau}_{max} = [1, ..., 1, 15, ..., 15]^T$. The reason is the same. $\boldsymbol{\tau}_{min}$ undercharges drivers while $\boldsymbol{\tau}_{max}$ overcharges drivers, both of which give rise to the highest objective function values and hence the lowest probability of finding the minimum solution in their proximity. To solve the minimization problem, the surrogate model naturally spends most of its effort exploring other regions in the design space, thereby predicting poorly for $\boldsymbol{\tau}_{min}$ and $\boldsymbol{\tau}_{max}$. Figure 5(b) shows the convergence history of the probabilistic EI. While exhibiting a brief increasing trend at the beginning of optimization, the pattern gradually and eventually converges to zero like Figure 3(b). Note that the probabilistic EI values in Figure 5(b) are generally smaller than those in Figure 3(b) because the probability of satisfying the constraint in Equation (30) is always less than or equal to one.

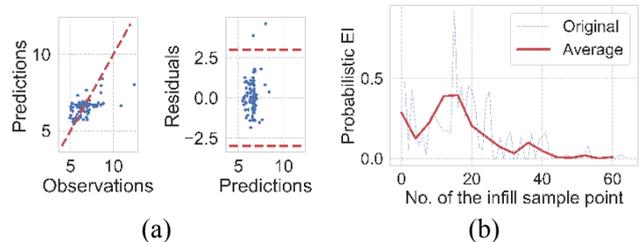

**Figure 5** Solving the bi-objective TLP: (a) validating the accuracy of the constructed surrogate model, and (b) convergence history of the probabilistic EI

Figure 6(a) shows the distribution of the 100 sample points based on their objective and constraint function values. Obviously, we are only interested in points lying below the constraint limit line represented by the blue filled circles. An interesting observation is that a Pareto front seems to appear suggesting a conflicting relation between the objective and the constraint. This observation, in part, supports our previous argument about the network flow drop. Specifically, while a higher toll price may decrease the objective function value, it may also increase the constraint function value by creating a more significant drop in the inflow to the PZ. A further reduced inflow equates to a more notable network flow drop and hence, a larger hysteresis loop in the NFD or a higher level of deviation from spread. Under the non-tolling scenario, the deviation from spread is the lowest as the PZ goes all the way to almost gridlock with no network unloading/recovery (see Figure 2(c)). The solution to the bi-objective TLP is shown in Figure 6(b) which corresponds to the corner point at the intersection of the Pareto front and the constraint limit line in Figure 6(a). The solution to the single-objective TLP is also shown by the green cross which has a lower objective function value but a higher constraint function value, as expected. Figure 6(c-f) show, respectively, the simulated averaged NFD, density, speed, and queue time series of the PZ under the optimal tolling scenario. The tolling NFD successfully maintains itself within the free-flow and at or near the capacity regimes without entering the congested branch of the non-tolling NFD. Traffic conditions in the PZ experience significant improvement with much lower densities and queues, and larger speeds.

### 4.4 Comparing the two optimal solutions

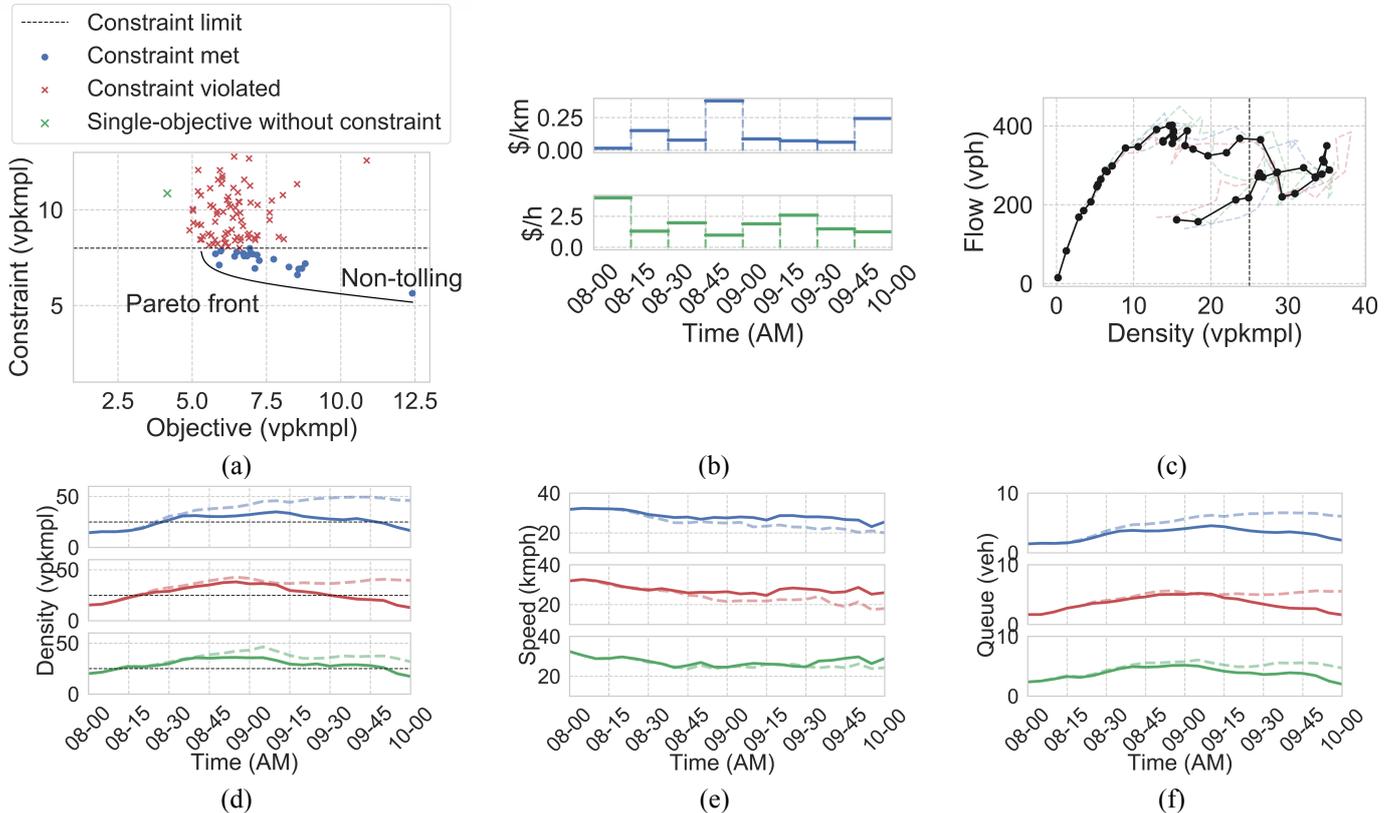

**Figure 6** (a) Distribution of the 100 sample points based on their objective and constraint function values, (b) optimal distance and delay toll patterns from solving the bi-objective TLP, (c) simulated NFDs for all replications (dashed lines) and the averaged NFD (solid line) of the PZ, (d-f) time series of density, speed, and queue of the PZ (dashed/solid lines represent the before-/after-pricing scenario)

Figure 7(a) compares the simulated NFDs of the PZ. The NFD from bi-objective optimization shifts more to the right because the heterogeneity constraint results in a lower toll price and hence a higher objective function value. Nevertheless, due to a lower constraint function value, higher flows are achieved during network loading which equates to a reduced network flow drop. During the transition period, although the NFD from bi-objective optimization works at higher densities, it produces similar or even slightly higher flows. Assuming a trapezoidal network exit function, there is a range of densities centering around the critical network density within which the flow can maintain at or near capacity (Daganzo, 2007; Mahmassani et al., 2013). Another observation is that the NFD from single-objective optimization exhibits a more significant local oscillatory loop. While the density remains almost constant, the flow undergoes a near-vertical jump along with a more heterogeneous distribution of congestion (see Figure 7(a) and (b)). This was also reported in Simoni et al. (2015). During network unloading/recovery, both NFDs exhibit a sizable hysteresis loop amplified by the very low demand entering the PZ at the end of simulation. Figure 7(b) shows that, although bi-objective optimization leads to higher densities, it produces slightly and consistently higher flows throughout the tolling period due to a lower level of the deviation from spread. Figure 7(c) and (d) show the average travel time in the PZ and in the entire network, respectively. Compared with the non-tolling scenario, the two optimal TLP solutions reduce the average travel time in the PZ by an average of 29.5% and 21.6%, respectively. Bi-objective optimization achieves less travel time improvement in the PZ because it allows the density to evolve further beyond the critical network density. While one may immediately question the 7.9% loss of travel time improvement in the PZ, a comparison between the average travel time in the entire network certainly provides the answer. Compared with the non-tolling scenario, bi-objective optimization reduces the average travel time in the entire network by an average of 2.5%, which is 1.1% higher than that by single-objective optimization. Hence, bi-objective optimization essentially manages to convert the 7.9% loss of travel time improvement in the PZ into the 1.1% gain of travel time improvement in the entire network. While producing less network-wide travel time improvement in the first two replications, single-objective optimization slightly increases the average network

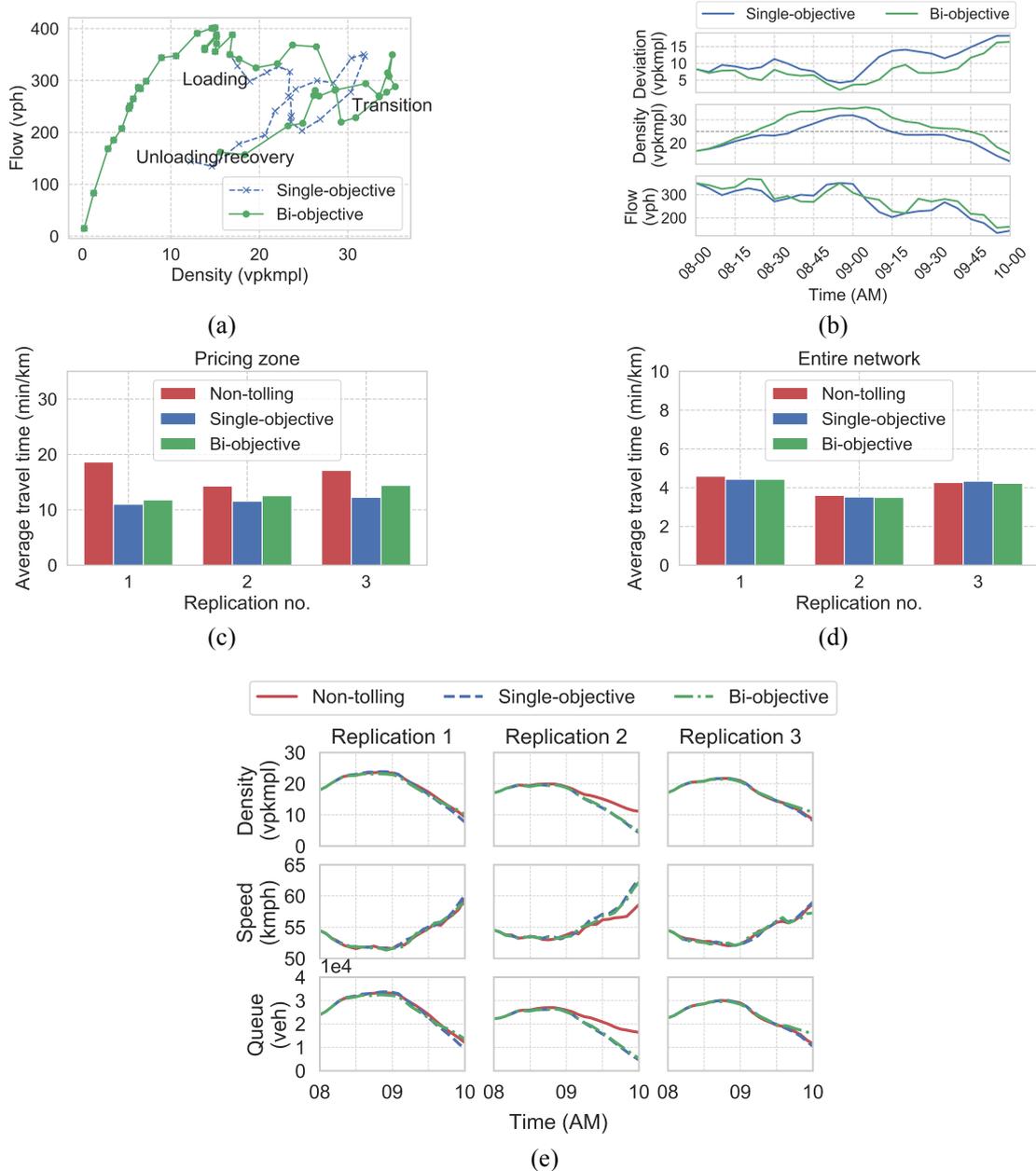

**Figure 7** Comparing the two optimal solutions: (a) simulated averaged NFDs of the PZ, (b) deviation, density, and flow time series of the PZ, (c) average travel time in the PZ, (d) average travel time in the entire network, and (e) density, speed, and queue time series of the entire network

travel time in the third replication probably due to overcharging the PZ and shifting congestion to the peripheral network. Two questions remain to be answered: (i) why is the travel time improvement in the entire network much lower than that in the PZ? And (ii) is it worthwhile to achieve the 1.1% gain of travel time improvement in the entire network at the cost of the 7.9% loss of travel time improvement in the PZ?

The answer to the first question is quite straightforward. The scale effect is a major reason given that the PZ only covers a relatively small area of the entire network (see Figure 2(a)). It is therefore no surprise that the performance of the entire network changes very little (see Figure 7(e)) when pricing a relatively small sub-network (Gu et al., 2018c). The performance of the entire network may even reduce, e.g. in the third replication, due to the redistribution of detour vehicles around the PZ which is highly dependent on the network configuration and the structure and magnitude of the demand, and hence case-specific. The proposed surrogate-based optimization framework represents an uncoordinated approach to pricing system design as our focus is explicitly and entirely on optimizing the performance of the PZ. Therefore, we need to check and ensure in an unsystematic manner that the optimal solution does not create unintended evident deterioration in the performance of the entire network. This is perhaps a limitation of our approach that motivates further investigation into two- or multi-area coordinated pricing.

The answer to the second question is a quick yes, at least from the authors' perspective. While acknowledging the fact that 1.1% is much lower and hence less seemingly appealing than 7.9%, we emphasize that the average travel time is normalized against the total distance traveled. Given that the total distance traveled in the entire network is over 60 times of that in the PZ, the total travel time saving in the entire network offered by the 1.1% is accordingly much higher than that in the PZ offered by the 7.9%. Indeed, we achieve on average a further network-wide total travel time saving of almost 700 hours. Hence from the entire network point of view, it is worthwhile to achieve the 1.1% gain of travel time improvement in the entire network at the cost of the 7.9% loss of travel time improvement in the PZ.

## 5 CONCLUSIONS

Toll optimization in a large-scale congested traffic network is often characterized by an expensive-to-evaluate objective function that requires considerable computational effort. An analytical solution or a brute force method is hardly possible as well as any commonly used metaheuristics. Therefore, in this paper, we solve two TLPs in a large-scale congested traffic network using surrogate-based optimization, a computationally efficient method for solving optimization problems that would otherwise involve much more expensive objective function evaluations. Specifically, we apply RK in conjunction with EI sampling as a global search method to determine the optimal toll pattern. Using the NFD to describe the level of network-wide congestion, we aim to optimize the time-varying JDDT such that the NFD of the PZ does not enter the congested regime. This logic formulates the first TLP which is a single-objective optimization problem. The second TLP is a direct extension by further considering minimizing the heterogeneity of congestion distribution. This is achieved by adding another objective in the formulation. To solve the bi-objective optimization problem, we convert one of the objectives into a constraint and reformulate a single-objective optimization problem.

A large-scale simulation-based DTA model of Melbourne, Australia is used to demonstrate the applicability of the surrogate-based toll optimization framework. Results show that both the optimal TLP solutions significantly improve the traffic conditions in the PZ including increased network speed, decreased network density, queue, and average travel time, and an NFD without the congested branch. A comparison further illuminates the advantage of bi-objective optimization. By considering and reducing the heterogeneity of congestion distribution, we achieve higher flows in the PZ as well as a lower average travel time or a higher total travel time saving in the entire network. We emphasize that surrogate-based optimization is also applicable to a wide range of high-dimensional network design problems other than toll optimization investigated in this paper.

The current research can be extended in a few promising directions. The first direction is to build and integrate a demand model with the DTA model to reflect more realistically drivers' behavioral changes in response to pricing. We do admit, however, that building a demand model itself is not a trivial task. The second direction, as has been touched upon early in the paper, is to investigate a coordinated approach that considers both the PZ and the peripheral network. This could be a major step towards designing a landmark coordinated multi-area pricing system. Note that the JDDT comprises a delay toll component and hence might raise safety concerns. An expedient is to apply speed limits and speeding penalties to discourage aggressive driving. To the best of our knowledge, different area-based pricing regimes have their respective pros and cons and there is no perfect one that currently exists in the literature. Further effort along this research line is therefore particularly desirable.

## ACKNOWLEDGMENTS

The authors would like to thank the Editor and the seven anonymous reviewers for their constructive comments and valuable suggestions to improve the quality of the paper.

## APPENDIX A

To investigate the effects of toll pattern smoothing parameters $\alpha$ and $\beta$ on the pricing control results, we perform a sensitivity analysis with two additional pairs of parameters: (i) $\alpha = \frac{1}{5}(1-0) = 0.2, \beta = \frac{1}{5}(15-0) = 3$, and (ii) $\alpha = \frac{1}{2}(1-0) = 0.5, \beta = \frac{1}{2}(15-0) = 7.5$. Theoretically speaking, a larger pair of $\alpha$ and $\beta$ imposes less constraint on the optimization and hence would achieve a better optimal objective function value, and vice versa. With a larger pair of $\alpha$ and $\beta$, the optimal toll patterns shown in Figure A.1(c) is, as expected, less smooth than those shown in Figure A.1(a). Accordingly, the simulated averaged NFD shown in Figure A.1(d) exhibits more chaotic behavior than that shown in Figure A.1(b) possibly due to the radical changes in the toll rates.

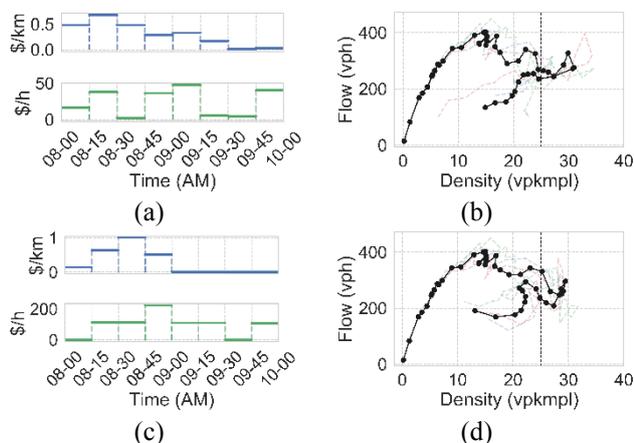

**Figure A. 1** Sensitivity analysis on the toll pattern smoothing parameters: (a) and (b) $\boldsymbol{\alpha = 0.2, \beta = 3}$, (c) and (d) $\boldsymbol{\alpha = 0.5, \beta = 7.5}$

## APPENDIX B

We further apply the surrogate-based optimization method to solve the single-objective TLP but with different percentages of adaptive driving in the simulation model including 60%, 80%, and 100%. Note that for this sensitivity analysis, we only consider the first replication of the simulation as it is far more time-consuming to run all the three replications. Also note that a higher percentage of adaptive driving generally leads to less congestion in the network and hence, we have to manually create some congestion by increasing the demand.

Figure B.1 clearly shows that the surrogate-based optimization method is robust and performs well for different assumed percentages of adaptive driving in the simulation model. The optimal objective function values as shown in Figure B.1(a) consistently reduce and converge as the number of function evaluations increases to 100. The before-and-after comparison of the simulated NFDs of the PZ as shown in Figure B.1(b-d) further confirms that the method works well across different adaptive driving scenarios and that the pricing control objective is consistently achieved (recall that 25 vpkmpl is the pricing control threshold).

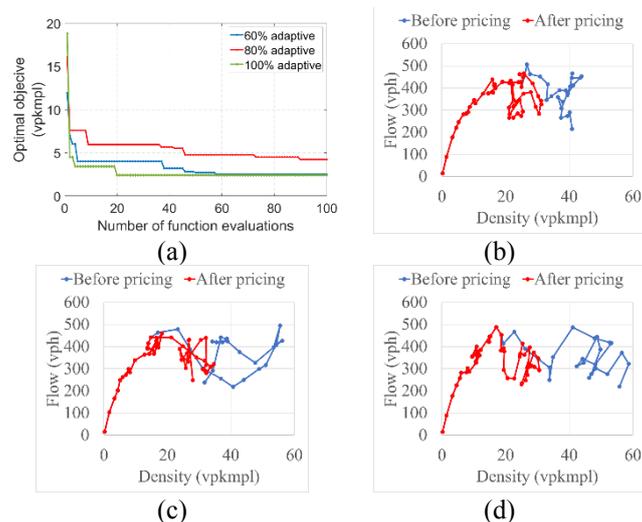

**Figure B. 1** Sensitivity analysis on the percentage of adaptive driving: (a) changes in the optimal objective function value as the number of function evaluations increases, and (b-d) 60%, 80%, and 100% adaptive driving